\documentclass{ifacconf}

\usepackage{natbib}

\usepackage{extract}
\usepackage{fancyhdr}
\usepackage{pdfpages}
\usepackage{amsmath}
\allowdisplaybreaks

\usepackage{bm}  
\usepackage{amsfonts}
\usepackage{amssymb}

\usepackage{mathtools} 
\usepackage{algorithm}
\usepackage{algpseudocode}
\usepackage{graphicx}
\usepackage{tabularx}
\usepackage{epsfig}
\usepackage{framed}
\usepackage{textcomp}
\usepackage{rotating}
\usepackage{booktabs}
\usepackage{pdflscape}
\usepackage[usenames,dvipsnames]{pstricks}
\usepackage{tikz}
\usepackage{enumerate}
\usetikzlibrary{mindmap,trees}

\usepackage{nomencl}
\usepackage{multicol}
\usepackage{verbatim}
\usepackage{etoolbox}

\renewcommand\nomgroup[1]{
  \item[\bfseries
  \ifstrequal{#1}{A}{Series}{%
  \ifstrequal{#1}{B}{Sets}{%
  \ifstrequal{#1}{C}{Other Symbols}{}}}%
]} 

\newtheorem{ass}{Assumption}

\definecolor{darkgreen}{RGB}{0,128,0}

\newcommand{\norm}[1]{\lVert{#1}\rVert}

\DeclareMathOperator*{\argmin}{arg\,min}

\DeclareMathOperator*{\minimise}{minimise \quad}

\DeclareMathOperator{\co}{Co}
\DeclareMathOperator{\vol}{vol}

\def\RR{\mathbb{R}}
\def\ZZ{\mathbb{Z}}

\def\PP{\mathbb{P}}
\def\EE{\mathbb{E}}

\def\N{\mathcal{N}}
\def\X{\mathcal{X}}
\def\U{\mathcal{U}}
\def\S{\mathcal{S}}

\def\W{\mathcal{W}}

\def\P{\mathcal{P}}
\def\probw{p_w}

\def\Id{{\mathcal {{I}}}}

\newcommand{\defeq}{:=}

\newcount\colveccount
\newcommand*\colvec[1]{
	\arraycolsep=1pt\def\arraystretch{1.5} 
        \global\colveccount#1
        \begin{pmatrix}
        \colvecnext
}

\def\colvecnext#1{
        #1
        \global\advance\colveccount-1
        \ifnum\colveccount>0
                \\
                \expandafter\colvecnext
        \else
                \end{pmatrix}
        \fi
}

\newtoks\rowvectoks
\newcommand{\rowvec}[2]{%
  \rowvectoks={#2}\count255=#1\relax
  \advance\count255 by -1
  \rowvecnexta}
  
\newcommand{\rowvecnexta}{%
  \ifnum\count255>0
    \expandafter\rowvecnextb
  \else
    \begin{pmatrix}\the\rowvectoks\end{pmatrix}
  \fi}
  
\newcommand\rowvecnextb[1]{%
    \rowvectoks=\expandafter{\the\rowvectoks&#1}%
    \advance\count255 by -1
    \rowvecnexta}

\newcommand{\bluetext}[1]{\textcolor{black}{#1}}

\allowdisplaybreaks[4]

\begin{document}
\begin{frontmatter}
\title{Adaptive robust predictive control with sample-based persistent excitation}

\author[Shanghai]{Xiaonan Lu}
and
\author[Oxford]{Mark Cannon}

\address[Shanghai]{Bosch China Investment Ltd, Shanghai, China}
\address[Oxford]{Department of Engineering Science, University of Oxford, UK\\
Email: \texttt{xiaonan.lu@cn.bosch.com}, \texttt{mark.cannon@eng.ox.ac.uk}}

\begin{abstract}
We propose a robust adaptive Model Predictive Control (MPC)
strategy with online set-based estimation for constrained
linear systems with unknown parameters and bounded
disturbances. A sample-based test applied to predicted
trajectories is used to ensure convergence of parameter
estimates by enforcing a persistence of excitation
condition on the closed loop system. The control law
robustly satisfies constraints and has guarantees of
feasibility and input-to-state stability. Convergence of
parameter set estimates to the actual system parameter
vector is guaranteed under conditions on reachability and
tightness of disturbance bounds.

  
\end{abstract}
\end{frontmatter}

\thispagestyle{empty}
\pagestyle{plain}

\section{Introduction}

Model Predictive Control (MPC) repeatedly solves a receding horizon optimisation problem to define a feedback control law. Controller performance relies on model accuracy, and mismatch between the controlled system and its model can degrade performance 
\citep{Mayne2000}. Various adaptive MPC algorithms have been proposed to improve model accuracy by estimating model parameters online. The control input of adaptive MPC has the dual purpose of regulating the controlled system and providing sufficient excitation to allow model identification.



Adaptive controllers rely on Persistent Excitation (PE) conditions to ensure asymptotic identification of model parameters~\citep{Narendra1987}, but combining PE conditions with state and control constraints is challenging~\citep{Mayne2014}. 
Several MPC algorithms have been developed with the aim of achieving both objectives. 
A common approach includes additional constraints in the MPC receding horizon optimisation problem. 
For example, \citet{Genceli1996} impose an additional periodic input constraint to guarantee a periodic persistently exciting feedback law for FIR systems, but the proposed algorithm has no guarantees of closed loop stability or convergence. 
\citet{Marafioti2014} consider ARMA models and formulate a PE constraint using past and current control inputs. 
The proposed algorithm ensures a PE condition if a periodic solution is found, but does not guarantee robust constraint satisfaction. 

An alternative approach uses controller performance measures to balance conflicting objectives of identifying model parameters and regulating the system state. 
\citet{Lu2020} propose an adaptive robust tube MPC formulation for linear systems that includes a PE measure in the predicted performance objective and suitably linearised PE conditions as constraints. This algorithm is proven to be recursively feasible and input-to-state stable, but it does not guarantee satisfaction of a PE condition in closed loop operation. 
 \citet{Heirung2015} avoid a nonconvex PE constraint by including an estimate of the parameter error covariance matrix in the MPC cost function. The optimisation problem is reformulated into a quadratically-constrained quadratic programming (QCQP) problem but recursive feasibility and stability were not proven.
Similarly, \citet{Weiss2014} include an information functional of the predicted parameter error covariance and enforce robust constraint satisfaction and stability by selecting the control input from the Robust Admissible Invariant (RAI) set, which results in a nonconvex optimisation problem.  
\bluetext{However, 
in both of these approaches there is no guarantee of persistently exciting control inputs.}

Other methods of promoting persistence of excitation in MPC have been proposed.
\citet{Bruggemann2022} reformulate the PE condition as an algebraic condition on a reference trajectory around a nominal equilibrium point. The algorithm ensures local convergence of the state and the parameter estimates. However, the approach does not \bluetext{guarantee robust  satisfaction of state and control constraints for all feasible initial conditions.}
\citet{Gonzalez2014} achieve stability and satisfaction of a PE condition using a switching control strategy. The algorithm drives the system state into a pre-defined invariant region for identification, within which an excitation signal is injected into the control law. The proposed algorithm is only applicable for open loop stable linear systems and the existence of the target region is example-dependent.

In this paper we ensure \citep[similarly to][]{Lu2022} persistence of excitation and recursive feasibility of the MPC optimisation. However, unlike \citet{Lu2022}, additional constraints to enforce persistent excitation are not included in the optimization because this would require nonconvex constraints on control inputs. Instead we propose a computationally tractable check on the solution by using postulated values (or samples) of
unknown model parameters and disturbances to generate sets of likely future trajectories and hence estimate the amount of excitation in past and future trajectories. For this check, persistence of excitation is evaluated over a single time-window rather than multiple windows. Each of these innovations reduces computation relative to~\citet{Lu2022}. 
%
For linear systems with uncertain parameters and bounded disturbances, the algorithm satisfies constraints robustly and achieves input-to-state \bluetext{practical} stability (ISpS) while retaining computational tractability.

The paper is organised as follows. 
Section~\ref{sec:adaptive_MPC} summarises the robust adaptive MPC formulation, explains set-based parameter estimation method, and states conditions for persistence of excitation. 
The sample-based check of a PE condition applied to predicted trajectories is described in Section \ref{sec:PE_sampled}, and Section \ref{sec:properties} shows that the control strategy ensures recursive feasibility, input-to-state \bluetext{practical} stability and asymptotic parameter convergence. 
Section~\ref{sec:example} illustrates the approach with a numerical example.

\textit{Notation:} $\mathbb{Z}_{\geq 0}$, $\mathbb{Z}_{> 0}$ are the sets of non-negative and positive integers, and $\mathbb{Z}_{[p,q]} = \{n\in\mathbb{Z} : {p \leq n \leq q}\}$.  
%
The Euclidean norm is $\norm{x}$ and $\|x\|_Q = (x^\top Q x)^{1/2}$.  
The identity matrix is $\Id$, and $A \succeq 0$ (or $A \succ 0$) indicates that $A$ is a positive semidefinite (positive definite) matrix.  
For $\mathcal{X},\mathcal{Y}\subset\mathbb{R}^n$, $A\mathcal{X} = \{Ax: x\in\mathcal{X}\}$, $\mathcal{X}\oplus\mathcal{Y} = \{x+y: x\in\mathcal{X},\, y\in\mathcal{Y}\}$ and $\co(\mathcal{X})$ is the convex hull of $\mathcal{X}$.
%
%
The value of a variable $y$ at the $k$th discrete time instant is denoted $y_k$ and $y_{k|t}$ denotes the value of $y_{k+t}$ predicted at time $t$.
Probabilities and expectations conditioned on the model state $x_t$ are denoted $\mathbb{P}(\cdot |x_t)=\mathbb{P}_t(\cdot)$ and $\EE(\cdot|x_t) = \EE_t(\cdot)$, and $\mathbb{P}(\cdot)$, $\EE(\cdot)$ are respectively equivalent to $\mathbb{P}_0(\cdot)$, $\EE_0(\cdot)$.

\section{Adaptive robust MPC}\label{sec:adaptive_MPC}


Consider systems described by the uncertain linear model
\begin{equation} \label{eq:sys}
x_{t+1} = A(\theta^\ast)x_t +B(\theta^\ast)u_t+ F w_t 
\end{equation}
with state, control input and unknown disturbance $x_t\in \RR^{n_x}$, $u_t \in \RR^{n_u}$ and $w_t \in\RR^{n_w}$ at discrete time $t\in\ZZ_{\geq 0}$. Here $\theta^\ast\in\mathbb{R}^p$ is an unknown parameter independent of $t$.

\begin{ass}[Disturbance] \label{ass:disturbance}
The disturbance sequence $\{w_t\in\W,\, t \in\ZZ_{\geq 0}\}$ is independent and identically distributed (i.i.d.), with $\EE (w_t) = 0$, $\EE (w_t w_t^\top) \succeq \epsilon_w \Id$, $\epsilon_w > 0$, and $\W$ is a known compact and convex set.
\end{ass}

\begin{ass}[Model parameters]\label{ass:param_uncertainty}
\begin{enumerate}[(a).]
\item $(A(\theta),B(\theta)) = (A_0,B_0) +\sum_{i = 1}^{p} (A_i,B_i) [\theta]_i$ for known $(A_i,B_i)$ \mbox{$i \in\ZZ_{[0,p]}$}.
\item $\Theta_0\ni \theta^\ast$ is a polytopic initial parameter set estimate.
\item $(A(\theta),B(\theta))$ is a reachable pair for all $\theta\in\Theta_0$.
\item $(A(\theta^\ast),B(\theta^\ast)) = (A(\theta),B(\theta))$ if and only if $\theta=\theta^\ast$. 
\end{enumerate}
\end{ass}

Starting from an initial set $\Theta_0$, 
we recursively refine a polytopic parameter set estimate $\Theta_t = \{ \theta : \Pi_{\Theta} \theta \leq \mu_t \}$ defined in terms of $\mu_t\in\mathbb{R}^{n_\Theta}$ and fixed $\Pi_{\Theta}\in\mathbb{R}^{n_{\Theta}\times p}$.

%

\subsection{Robust MPC}
Consider a predicted control law parameterised at time $t$ in terms of decision variables $\mathbf{v}_t = \{v_{0|t},\ldots,v_{N-1|t}\}$ as
\begin{align} \label{eq:input_law}
u_{k|t} =\begin{cases}
Kx_{k|t}+v_{k|t}, & k \in \ZZ_{[0,N-1]}\\
Kx_{k|t}+s_{k|t}, & k \geq N
\end{cases}
\end{align}
where $N$ is the MPC prediction horizon and the gain $K$ is such that $z_{t+1} \in \co\{A_K(\theta) z_t,\, \theta\in\Theta_0\} $ is quadratically stable with $A_K(\theta)=A(\theta) + B(\theta)K$. Here $\{s_{k|t},\, k\in \ZZ_{\geq 0}\}$ is a zero-mean i.i.d.\ random sequence that is introduced into the terminal control law to ensure persistent excitation, with $s_{k|t}\in\S$ for some known compact and convex set $\S$, and  $\EE(s_{k|t}\smash{s_{k|t}^\top}) \succeq \epsilon_s\Id \succ 0$ for all $k,t\in\ZZ_{\geq 0}$.

%
%
We assume linear state and control input constraints:
\begin{equation} \label{eq:input_state_constraint}
  x_{k|t} \in \mathcal{X}, \quad u_{k|t} \in \U , \quad 
\forall k \geq 0 ,
\end{equation}
where $\X$, $\U$ are given polytopic sets.
To enforce these constraints we define a terminal set $\X_T$ satisfying 
\begin{align}
&\X_T \subseteq\X, \quad K\X_T \oplus \S \subseteq \U ,
\label{eq:terminala}
\\
&A_K(\theta) \X  \oplus B(\theta)\S \oplus F \W \subseteq \X_T \quad
\forall \theta \in \Theta_0
\label{eq:terminalb}
\end{align}
%
%
To ensure satisfaction of (\ref{eq:input_state_constraint})
we construct a tube $\mathbf{X}_t$,
\[
\mathbf{X}_t = \{\X_{k|t},\, k\in\mathbb{Z}_{[0,N]}\},
\]
where the sets $\X_{k|t}$ satisfy, for all $\theta\in\Theta_t$, $k\in\mathbb{Z}_{[0,N-1]}$,
\begin{align}
& \X_{k|t} \subseteq \X, \quad K \X_{k|t} \oplus \{ v_{k|t}\} \subseteq \U ,
\label{eq:tube_xu_constraint}
\\
& A_K(\theta)\X_{k|t} \oplus \{ B(\theta)v_{k|t}\} \oplus F\W \subseteq \X_{k+1|t} ,
\label{eq:recurrence_constraint}
\end{align}
and the initial and terminal conditions 
\begin{equation}\label{eq:terminal_constraint}
x_t \in \X_{0|t}, \quad
\X_{N|t} \subseteq \X_T .
\end{equation}
%

\begin{rem}\label{rem:ell_tubes}
The sets comprising $\mathbf{X}_t$ may be defined, for example, as homothetic or fixed-complexity polytopic sets~\citep{Lorenzen2018, Lu2019}. However, in the example of Section~\ref{sec:example}, we use ellipsoidal sets $\X_{k|t} = \{x\in\smash{\mathbb{R}^{n_x}}: \smash{(x-z_{k|t})^\top P (x-z_{k|t})} \leq \beta_{k|t}\}$, which ensure favourable scaling of computation with $n_x$ when $P$ is computed offline and $\boldsymbol{\beta}_t = \{\beta_{0|t},\ldots, \beta_{N|t}\}$ and ${\bf z}_t = \{z_{0|t},\ldots,z_{N|t}\}$ are online decision variables.
\end{rem}

%

We consider a  nominal performance objective defined for a given parameter vector $\smash{\hat{\theta}_t}\in\Theta_t$ by the cost function
\[
 J(x_t, \hat{\theta}_t, {\bf v}_t) = \sum_{k = 0}^{N-1} l(\hat{x}_{k|t},K\hat{x}_{k|t} + v_{k|t}) + V_{N|t}(\hat{x}_{N|t},\hat{\theta}_t)
%
\]
with $\smash{\hat{x}_{0|t}} = x_t$, $\smash{\hat{x}_{k+1|t}} = \smash{A_K(\hat{\theta}_t) \hat{x}_{k|t} + B(\hat{\theta}_t) v_{k|t}}$, $k\in\mathbb{Z}_{[0,N-1]}$. The stage cost $l(\cdot,\cdot)$ and terminal cost $V_{N|t}(\cdot,\smash{\hat{\theta}_t})$ are assumed to be positive definite quadratic functions
%
satisfying, for given $\smash{\hat{\theta}_t}\in\Theta_t$ and all $x\in\mathbb{R}^{n_x}$,
\begin{equation}\label{eq:terminal_cost}
{V_{N|t}(x,\hat{\theta}_t)  = V_{N|t}\bigl(A_K(\hat{\theta}_t)x,\hat{\theta}_t \bigr) +l(x, Kx).}
\end{equation}
The nominal parameter vector $\smash{\hat{\theta}_t}\in\Theta_t$ may be obtained by projecting a least squares parameter estimate onto $\Theta_t$ \citep{Lorenzen2018} or by projecting $\smash{\hat{\theta}_{t-1}}$ onto $\Theta_t$:
\begin{equation} \label{eq:nominal_theta_update}
\hat{\theta}_{t} = \argmin_{\theta\in\Theta_t} \|\hat{\theta}_{t-1} - \theta\| .
\end{equation}
A robust tube MPC law is determined by the solution, denoted $(\mathbf{v}_t^o, \mathbf{X}_t^o)$, of the problem of minimizing $J(x_t,\hat{\theta}_t,\mathbf{v}_t)$ over $\mathbf{v}_t$ and $\mathbf{X}_t$ subject to  (\ref{eq:tube_xu_constraint})-(\ref{eq:terminal_constraint}). This is the basis of the algorithm defined in Section~\ref{sec:PE_sampled}.

%

\subsection{Set-based parameter estimate}\label{sec:parameter_set_estimation}


%

Set membership identification uses information on $u_{t-1}$, $x_t$ and $x_{t-1}$ to define a set $\Delta_t$ of unfalsified model parameters at time $t$. This is combined with the parameter set estimate $\Theta_{t-1}$ to construct an updated estimate $\Theta_t$ as follows. 
The linear dependence of the system model (\ref{eq:sys}) on $\theta^\ast$ implies
\[
x_{t+1} = \Phi(x_t, u_t) \theta^\ast + \phi(x_t, u_t)+ F w_t ,
\]
where $\Phi_t=\Phi(x_t, u_t)$ and $\phi_t = \phi(x_t, u_t)$ are defined 
\begin{align}
\Phi_t & =\Phi(x_t,u_t)  = \begin{bmatrix} A_1 x_t\!+\!B_1 u_t &  \cdots & A_p x_t\! + \! B_p u_t\end{bmatrix} 
\label{eq:Phik}
\\
\phi_t & =\phi(x_t,u_t) = A_0 x_{t}+B_0 u_{t} . 
\end{align}
A set of unfalsified parameters at time $t$ given $x_t$, $x_{t-1}$, $u_{t-1}$, and the disturbance set $\W$, is given by
\[
\Delta_t = \{ \theta : x_t - \Phi_{t-1}\theta - \phi_{t-1} \in  F\W\} .
\]
Using $N_\mu$ unfalsified sets, the minimum volume $\Theta_t$ satisfying $\Theta_t\supseteq \smash{\bigcap_{\tau = t-N_\mu+1}^t}\Delta_\tau \cap \Theta_{t-1}$ is found by solving $n_{\Theta}$ convex programs to compute the update $\mu_t$. This ensures that $\theta^\ast\in\Theta_t\subseteq \Theta_{t-1}$ for all~$t\in\mathbb{Z}_{> 0}$~\citep{Lu2020}.


%
\subsection{Persistent excitation}\label{sec:PE_linear_feedback}

The regressor $\Phi_t$ in (\ref{eq:Phik}) is persistently exciting \citep{Narendra1987} if there exists a horizon $N_u$ and scalar $\epsilon_\Phi >0$ such that
$\sum_{k=t}^{t+N_u-1} \Phi_k^\top \Phi_k  \succeq \epsilon_\Phi \Id$ for all $t \in \ZZ_{\geq 0}$.
In the current work, we define persistent excitation using the requirement that there exists $\epsilon_\Phi > 0$ and an infinite sequence of time instants $t_0,t_1,\ldots$ such that
\begin{equation}\label{eq:pe_condition}
\sum_{k=t_i}^{t_i+N_u-1} \Phi_k^\top \Phi_k \succeq \epsilon_\Phi \Id 
\quad \forall \, i\in \ZZ_{\geq 0} .
\end{equation}
We refer to the interval $\mathbb{Z}_{[t_i, t_i+N_u-1]}$ as a \textit{PE window}.

\begin{ass}[Tight disturbance bound]\label{ass:probw}
For all $w^0\in \partial\W$ and any $\epsilon > 0$ the disturbance sequence $\{w_0,w_1,\ldots\}$ satisfies $\PP \bigl\{ \norm{ w_t - w^0 } < \epsilon \bigr\} \geq \probw (\epsilon) $, for all $t\in\ZZ_{\geq 0}$, where $\probw(\epsilon) > 0$ whenever $\epsilon > 0$.
\end{ass}
The PE condition~(\ref{eq:pe_condition}) ensures convergence of the parameter set $\Theta_t$ 
(see e.g.~\citet{Lu2022} Lemma~1).

\begin{lem}\label{lem:convergence}
Under Assumptions~\ref{ass:disturbance} and \ref{ass:probw},
if~(\ref{eq:pe_condition}) holds and $ N_u \leq N_\mu $, then
$\smash{\displaystyle\lim_{t\to\infty} \Theta_t} = \{\theta^\ast\}$ with probability 1 (w.p.\,1).
\end{lem}

\section{Persistence of excitation check}\label{sec:PE_sampled}

The PE condition (\ref{eq:pe_condition}) is not suitable to be included directly as a constraint in the optimization of predicted trajectories because it involves nonconvex quadratic matrix inequalities. Instead, we impose (\ref{eq:pe_condition})  with a certain level of statistical confidence using a sample-based check applied to the solution of the MPC optimization. This approach performs a comparison with a reference solution that is known to satisfy the PE condition on average, and adopts the reference as a fallback solution if the check fails. 

Given the solution, $(\mathbf{v}_t^o, \mathbf{X}_t^o)$, of the MPC optimization at time $t$, we compare a measure of PE for $\mathbf{v}^o_t$ and a reference solution $\hat{\mathbf{v}}_t$ over a set of samples of model parameters $\theta_t\in\Theta_t$ and disturbance sequences $\mathbf{w}_t=\{w_{0|t},\ldots,w_{N-1|t}\}\in\smash{\W^{N}}$.  The reference $\hat{\mathbf v}_t=\{\smash{\hat{v}_{0|t}},\ldots,\smash{\hat{v}_{N-1|t}}\}$, is defined using the solution  $\mathbf{v}^o_{t-1}$ 
adopted at time $t-1$ and the random input $s_{N-1|t}$: 
\[
\hat{v}_{k|t} = \begin{cases}
    v^o_{k+1|t-1} ,  & k\in\mathbb{Z}_{[0,N-2]} \\
    s_{N-1|t} , & k=N-1.  \end{cases} 
\]

Let $\smash{\zeta_t^{(i)}}=(\smash{\theta_t^{(i)}}, \smash{ w_{0|t}^{(i)}},\ldots, \smash{ w_{N-1|t}^{(i)}})$
be a sample from the uniform distribution on $\Theta_t\times \W^{N}$, and let 
$\{\smash{\zeta_t^{(1)}},\ldots, \smash{\zeta_t^{(N_s)}}\}$\rule{0pt}{9.5pt}
%
%
denote a set of $N_s$ independent samples. Corresponding predicted state and control sequences $(\smash{\tilde{x}_{k|t}^{(i)}},\smash{\tilde{u}_{k|t}^{(i)}})$, \bluetext{for $k \in \ZZ_{[0, N-1]}$, are generated using}
\begin{equation}\label{eq:xu_def}
\begin{aligned}
\tilde{x}_{k+1|t}^{(i)} &= A(\theta^{(i)}_t) \tilde{x}_{k|t}^{(i)} + B(\theta^{(i)}_t) \tilde{u}_{k|t}^{(i)} + F w^{(i)}_{k|t} ,
\\
\tilde{u}^{(i)}_{k|t} & = K\tilde{x}_{k|t}^{(i)} + v_{k|t} ,
\end{aligned}
\end{equation}
and we set $\tilde{x}_{k|t}^{(i)} = x_{k+t}$ if $k \leq 0$ and $\tilde{u}_{k|t}^{(i)} = u_{k+t}$ if $k<0$.
For a PE window starting at $t+\kappa$ with $\kappa\in\ZZ_{[-N_u+1,N-N_u]}$, the associated \textit{PE matrix} $\Psi_{\kappa|t}(\mathbf{v}_t, \smash{\theta^{(i)}_t}, \smash{\mathbf{w}^{(j)}_t})$ is defined as
\[
\Psi_{\kappa|t}(\mathbf{v}_t, \zeta^{(i)}_t) = \sum_{k=\kappa}^{\kappa+N_u-1}
\Phi(\tilde{x}_{k|t}^{(i)},\tilde{u}_{k|t}^{(i)})^\top \Phi(\tilde{x}_{k|t}^{(i)},\tilde{u}_{k|t}^{(i)}) .
\]
We also allow the PE window to extend beyond the first $N$ time steps of the prediction horizon by defining $v_{k|t} = s_{k|t}$ and $\smash{w_{k|t}^{(i)} = w_{k+t}}$ for all $i\in\ZZ_{[1,N_s]}$ and all $k\geq N$ in (\ref{eq:xu_def}), and defining
\begin{align*}
\Psi_{\kappa|t}(\mathbf{v}_t, \zeta^{(i)}_t) &= \sum_{k=\kappa}^{N-1}
\Phi(\tilde{x}_{k|t}^{(i)},\tilde{u}_{k|t}^{(i)})^\top \Phi(\tilde{x}_{k|t}^{(i)},\tilde{u}_{k|t}^{(i)})  
\\
&\quad +
\sum_{k = N}^{\kappa + N_u - 1} \EE_{N|t} \bigl[\Phi(\tilde{x}_{k|t}^{(i)},\tilde{u}_{k|t}^{(i)})^\top \Phi(\tilde{x}_{k|t}^{(i)},\tilde{u}_{k|t}^{(i)})  \bigr]
\end{align*}
for $\kappa \in \ZZ_{[N-N_u+1,N-1]}$. Here $\EE_{N|t} (\cdot )$ 
denotes the expectation (conditioned on $\smash{\tilde{x}_{N|t}}$)  that is obtained  by marginalising over the distributions of
$s_{k|t}$ and $w_{k+t}$ for $k\geq N$.

To perform the check, we compare 
$\delta_{\kappa|t}$ and $\smash{\hat{\delta}_{\kappa|t}}$ defined by
\begin{equation}\label{eq:delta}
\begin{aligned}[t]
\delta_{\kappa|t} & = \max_{\delta\in\mathbb{R}} \, \delta \ \text{subject to} 
\ 
\Psi_{\kappa|t}(\mathbf{v}^o_t, \zeta_t^{(i)}) \succeq \delta\Id 
\ \forall 
i \in \ZZ_{[1,N_{s}]}
\\
\hat{\delta}_{\kappa|t} & = \max_{\delta\in\mathbb{R}} \, \delta \ \text{subject to} 
\ 
\Psi_{\kappa|t}(\hat{\mathbf{v}}_t, \zeta_t^{(i)}) \succeq \delta\Id 
\ \forall 
i \in \ZZ_{[1,N_{s}]}
\end{aligned}
\end{equation}
The check is said to have failed if $\delta_{\kappa|t} < \hat{\delta}_{\kappa|t}$ and in this case we redefine $\mathbf{v}^o_t$ by setting its value to $\hat{\mathbf{v}}_t$.

At each time $t$ we propose to check a single PE window with start time $t+\kappa$, and to decrease $\kappa$ by $1$ at successive time steps until $\kappa = -N_u+1$,  when it is reset to $\kappa=  N-1$. The overall strategy of Algorithm~\ref{alg:sampling} is therefore to track a PE window with a fixed start time as time advances until the window contains only the past and current time instants from $t-N_u+1$ to $t$. We show in Sections~\ref{sec:properties} and~\ref{sec:example} that this strategy ensures that the closed loop system has a non-zero probability of satisfying the PE condition~(\ref{eq:pe_condition}) at least once every $N+N_u-1$ time steps, and therefore~(\ref{eq:pe_condition}) necessarily holds at an infinite number of time instants.

\begin{algorithm}
\caption{Adaptive MPC with sampled PE check}\label{alg:sampling}
\textbf{Offline:}
Compute $\X_T$ satisfying (\ref{eq:terminala})-(\ref{eq:terminalb}) and set $\kappa = N-1$.\\
\textbf{Online:} At times $t=0,1,\ldots$:
\begin{enumerate}[(a).]
\item Obtain the current state $x_t$.
\item Update $\Theta_{t} $ and $\smash{\hat{\theta}}_t$, and determine $V_{N|t}$ satisfying (\ref{eq:terminal_cost}).
\item 
Compute the solution $(\mathbf{v}^o_t, \mathbf{X}_t^o)$ of the convex program:
\[
\begin{aligned}
\P(x_t,\Theta_t,\hat{\theta}_t):\   \minimise_{\mathbf{v_t},\mathbf{X}_t} \ 
J(x_t,\smash{\hat{\theta}}_t,\mathbf{v}_t) \ \text{s.t.\ (\ref{eq:tube_xu_constraint})-(\ref{eq:terminal_constraint})}.
\end{aligned}
\]
\item 
Draw samples $\{\smash{\zeta^{(1)}},\ldots,\smash{\zeta^{(N_s)}}\}$
and compute $\delta_{\kappa|t}$, $\smash{\hat{\delta}_{\kappa|t}}$.\\
Set $\mathbf{v}_t^o \gets \smash{\hat{\mathbf{v}}_t}$ if $\smash{\delta_{\kappa|t}} < \smash{\hat{\delta}_{\kappa|t}} + \varepsilon$ for a given threshold $\varepsilon$.
\item
Apply the control input $u_t = K x_t + v_{0|t}^o$.
\\
If $\kappa > -N_u+1$ set $\kappa \gets \kappa - 1$, otherwise set $\kappa \gets N-1$.
\end{enumerate}
\end{algorithm}

\begin{rem}
The threshold $\epsilon$ in step (d) is a positive value used to distinguish non-zero minimum eigenvalues in (\ref{eq:delta}) from those that are zero to numerical precision. 
\end{rem}

\begin{rem}
The computational effort to find $\smash{\delta_{\kappa|t}}$ and $\smash{\hat{\delta}_{\kappa|t}}$ in (\ref{eq:delta}) is equivalent to  determining the minimum eigenvalues of $2N_{s}$ positive semidefinite $p\times p$ matrices.
\end{rem}

\begin{rem}
If the disturbance bound $\W$ is ellipsoidal and the sets $\X_{k|t}$ are ellipsoidal (as described in Remark~\ref{rem:ell_tubes}), then the update of $\Theta_t$ in step (b) and problem $\P$ in step (c) become second-order conic programs (SOCPs).
\end{rem}

\section{Closed loop system properties}\label{sec:properties}


The MPC optimisation (problem $\P$ in Algorithm \ref{alg:sampling})
is by construction feasible at all times $t > 0$ if it is feasible at $t = 0$ since $\Theta_{t}\subseteq\Theta_{t-1}$ for all $t$.
In addition, the closed loop system has the following input-to-state stability property (see e.g.~Theorem 1 and Corollary 1 of~\citet{Lu2020}).


\begin{prop}[Stability]\label{prop:stability}
Under Assumptions \ref{ass:disturbance} and \ref{ass:param_uncertainty}, the system (\ref{eq:sys})  with Algorithm~\ref{alg:sampling} is input-to-state practically stable (ISpS)~\citep[][Def.~6]{Limon2008}
with respect to an input defined by $Fw_t + B(\theta^\ast)s_t$,
where $s_t = s_{N|t-N}$ for $t \geq N$ and $s_t = 0$ for $t<N$,
in the set of initial conditions $x_0$ for which $\mathcal{P}(x_0,\Theta_0,\smash{\hat{\theta}_0})$ is feasible.
\end{prop}


A consequence of Proposition~\ref{prop:stability} \citep{Limon2008,Lu2020} is that a $\mathcal{K}\mathcal{L}$-function $\eta$ and $\mathcal{K}$-functions $\psi$, $\xi$ exist such that, 
\bluetext{for all $t\in\ZZ_{\geq 0}$ and all feasible $x_0$},
\begin{multline*}
\norm{x_t} \leq \eta(\|x_0\|, t) + \psi\bigl(\max_{\tau\in\ZZ_{[0,t-1]}} \|Fw_\tau + B(\theta^\ast)\bluetext{s_{\tau}}\| \bigr)
\\
+ \xi\bigl( \max_{\tau\in\ZZ_{[0,t-1]}} \| \smash{\hat{\theta}_t} - \theta^\ast \|\bigr) .
\end{multline*}


\subsection{Closed loop system persistence of excitation} \label{sec:properties2}

This section considers the randomised check performed in step (d) of Algorithm~\ref{alg:sampling}. We show that this procedure ensures that the closed loop system satisfies the PE condition (\ref{eq:pe_condition})
%
using an approach that does not require knowledge of the probability distribution of the disturbance $w_t$ or the prior distribution of the unknown model parameters $\theta^\ast$.
To simplify notation, let $\smash{d_{\kappa|t}(\mathbf{v})}$ denote the minimum over $\zeta\in\Theta_t\times \W^{N}$ of the smallest eigenvalue of $\Psi_{\kappa|t}(\mathbf{v},\zeta)$.

\begin{lem}\label{lem:delta}
If $\smash{\delta_{\kappa|t}} \geq \smash{\hat{\delta}_{\kappa|t}} + \varepsilon$, then, for some $\rho \in [0,1)$,
\begin{equation}\label{eq:expp_delta}
\mathbb{P}_t \bigl\{ d_{\kappa|t} (\mathbf{v}_t^o) \geq d_{\kappa|t} (\hat{\mathbf{v}}_t)\bigr\} \geq 1-\rho^{N_s} .
\end{equation}
\end{lem}

\textbf{Proof.}\hspace{1ex}%
The minimum eigenvalue of $\smash{\Psi_{\kappa|t}}(\mathbf{v}_t^o,\zeta)$ is a Lipschitz continuous function of $\zeta$ since $\{(A_i,B_i), \, i\in\ZZ_{[1,p]}\}$, $\X$, $\U$, $\Theta_0$ and $N_{u}$ are by assumption bounded. 
Therefore $\smash{\delta_{\kappa|t}} \geq \smash{\hat{\delta}_{\kappa|t}} + \varepsilon$ implies $ d_{\kappa|t} (\mathbf{v}_t^o) \geq \smash{\hat{\delta}_{\kappa|t}}$ whenever $\smash{\zeta^{(i)}_t}\rule{0pt}{9.3pt}$ for some $i\in \ZZ_{[1,N_s]}$ lies in a certain neighbourhood, $\N$, of a minimizing argument $\smash{\zeta^\ast_t}$ such that $d_{\kappa|t}(\mathbf{v}_t^o) = \Psi_{\kappa|t}(\smash{\mathbf{v}_t^o,\zeta^\ast_t})$
(where $\N$ depends on $\{(A_i,B_i), \, i\in\smash{\ZZ_{[1,p]}}\}$, $\X$, $\U$, $\Theta_0$, $N_{u}$, and $\varepsilon$, \bluetext{and $\N$ necessarily has a non-empty interior}). 
But the samples $\smash{\zeta^{(i)}_t}\rule{0pt}{9.3pt}$, $i\in\smash{\ZZ_{[1,N_s]}}$ are uniformly 
distributed over $\Theta_t\times\W^N$.
Denoting the probability that $\smash{\zeta^{(i)}_t}\rule{0pt}{9.3pt}$ lies outside $\N$ as $\rho$, where $\rule{0pt}{9.3pt}\rho \leq 1-\vol\{\N\cap(\Theta_t\times\W^N)\}/\vol\{\Theta_t\times\W^N\} \in [0,1)$, 
it follows that $\smash{\zeta^{(i)}_t}\in \N \rule{0pt}{9.3pt}$ for one or more index $i\in\ZZ_{[1,N_s]}$ with probability $1-\rho^{N_s}$
%
Therefore (\ref{eq:expp_delta}) must hold because $\smash{\hat{\delta}_{\kappa|t}} \geq d_{\kappa|t}(\hat{\mathbf{v}}_t)$.
\qed

The following result demonstrates that the PE condition (\ref{eq:pe_condition}) holds under Algorithm~\ref{alg:sampling}, and therefore that $\Theta_t$ converges to $\{\theta^\ast\}$ with probability 1 by Lemma~\ref{lem:convergence}.

\begin{thm}\label{thm:pe_closed_loop}
Suppose Assumptions \ref{ass:disturbance} and \ref{ass:param_uncertainty} are satisfied and $N_u > n_x$,
then Algorithm~\ref{alg:sampling} ensures that (\ref{eq:pe_condition}) holds at an infinite number of instants $t_0,t_1,\ldots$, for some $\epsilon_\Phi>0$.
\end{thm}

\textbf{Proof.}\hspace{1ex}%
For all $\theta \in\RR^p$ and any $\zeta^{(i)}\in\Theta_t \times \W^N$ we have  
\[
\EE_{N|t} \bigl[ \theta^\top \Psi_{N-1|t} (\smash{\hat{\mathbf{v}}}_{t},\zeta^{(i)}) \, \theta \bigr]
=
\!\!\! \sum_{k=N-1}^{N+N_u-2} \!\!\!\! \EE_{N|t}\bigl[ \|\Phi(\tilde{x}^{(i)}_{k|t},\tilde{u}^{(i)}_{k|t}) \theta\|^2\bigr],
\]
where $\smash{\tilde{x}_{k|t}^{(i)}}$ and $\smash{\tilde{u}_{k|t}^{(i)}}$ are generated by (\ref{eq:xu_def}) with $v_{k|t} = s_{k|t}$.
Therefore, from Assumptions~\ref{ass:param_uncertainty}(c) and \ref{ass:param_uncertainty}(d), and $\Theta_t\subseteq\Theta_0$, $\epsilon_s > 0$, and $\hat{v}_{N-1|t}=s_{N-1|t}$, it follows that $\epsilon_\Phi > 0$ exists so that (e.g.\ see \citet{Lu2022}, Thm.~3):
\begin{equation}\label{eq:epsilon_Phi}
\EE_{N|t} \bigl[ \theta^\top\Psi_{N-1|t} (\hat{\mathbf{v}}_t,\zeta^{(i)}) \, \theta \bigr]
\geq \epsilon_\Phi\|\theta\|^2 .
\end{equation}
Hence the probability that the minimum eigenvalue of $\Psi_{N-1|t} (\smash{\hat{\mathbf{v}}_t},\zeta)$ for all $\zeta\in\Theta_t \times \W^N$ is no less than $\epsilon_\Phi$ satisfies
\begin{equation}\label{eq:pe_closed_loop_cond1}
\PP_{t} \bigl[ d_{N-1|t} (\smash{\hat{\mathbf{v}}_t}) \geq \epsilon_\Phi \bigr] \geq \alpha 
\end{equation}
for some $\alpha \in (0,1]$. 
This property of predicted trajectories ensures the PE condition~(\ref{eq:pe_condition}) because
Algorithm~\ref{alg:sampling} implies
\begin{equation}
\PP_{t} \bigl[ d_{\kappa|t} (\mathbf{v}_t^o) 
\geq
d_{\kappa|t} (\hat{\mathbf{v}}_t)  \bigr]
\geq 1-\rho^{N_s} \in (0,1]
\end{equation}
(this follows from Lemma~\ref{lem:delta} if $\delta_{\kappa|t} \geq \smash{\hat{\delta}_{\kappa|t}} + \epsilon$ and from $\mathbf{v}^o_t=\hat{\mathbf{v}}_t$ if $\delta_{\kappa|t} < \smash{\hat{\delta}_{\kappa|t}} + \epsilon$), and hence, for any $\beta \in (0,1]$,
\begin{equation}\label{eq:pe_closed_loop_cond2}
\PP_t \bigl[ d_{\kappa|t} (\hat{\mathbf{v}}_{t}) \geq
\epsilon_\Phi\bigr]  \geq \beta
\ \Rightarrow \ 
\PP_t \bigl[ d_{\kappa|t} (\mathbf{v}_t^o) \geq \epsilon_\Phi \bigr] \geq \beta (1-\rho^{N_s}) .
\end{equation}
Furthermore, at time $t$, 
$\Psi_{\kappa-1|t+1} (\hat{\mathbf{v}}_{t+1},\zeta)$ 
and $\Psi_{\kappa|t} (\mathbf{v}_{t}^o,\zeta)$ 
are identically distributed random variables for all $\kappa\in\ZZ_{[-N_u+1, N-1]}$, so that, for any $\beta \in (0,1]$,
\begin{equation}\label{eq:pe_closed_loop_cond3}
\PP_t \bigl[ d_{\kappa|t} (\mathbf{v}_{t}^o ) \geq
\epsilon_\Phi \bigr]  \geq \beta
\ \Rightarrow\  
\PP_t \bigl[ d_{\kappa-1|t+1} (\hat{\mathbf{v}}_{t+1}) \geq \epsilon_\Phi \bigr] \geq \beta.
\end{equation}
Let $N_{pe} := N + N_u - 1$, starting from (\ref{eq:pe_closed_loop_cond1}) and applying (\ref{eq:pe_closed_loop_cond2})-(\ref{eq:pe_closed_loop_cond3}) 
$N_{pe}$
times, we obtain
\[
\PP_t \bigl[ d_{-N_u+1|t+N_{pe} } (\mathbf{v}_{t+N_{pe}}^o,\zeta ) \geq \epsilon_\Phi \bigr]
\geq (1-\rho^{N_s})^{N_{pe}}\alpha .
\]
But $\theta^\top\Psi_{-N_u+1|t+N_{pe}} (\mathbf{v}_t^o,\zeta)\, \theta = \sum_{k=N-1}^{N+N_u-2} \|\Phi_{t+k}\theta\|^2$ for all $\theta\in\RR^{p}$, where $\Phi_{t+k} = \Phi(x_{t+k},u_{t+k})$ is the regressor defined in (\ref{eq:Phik}) evaluated along the closed loop trajectories of (\ref{eq:sys}) under Algorithm~\ref{alg:sampling}.
Therefore, for all $\theta\in\RR^{p}$,
\[
\PP_t \Bigl( \sum_{k=N-1}^{N+N_u-2} \|\Phi_{t+k} \theta\|^2\geq \epsilon_\Phi \|\theta\|^2\Bigr)
\geq (1-\rho^{N_s})^{N_{pe}}\alpha > 0 .
\]
Since this holds at $t = nN_{pe}$, $n\in\ZZ_{>0}$, the PE condition (\ref{eq:pe_condition}) must be satisfied at an infinite number of instants.
\qed

\section{Numerical examples} \label{sec:example}

We consider a system with 3 unknown parameters, defined by the following matrices in (\ref{eq:sys}) and Assumption~\ref{ass:param_uncertainty}(a):
\begin{align*}
\hspace{-0.1ex} (A_0, B_0, F)  &= \biggl[\!\begin{smallmatrix*}[r] 
1.08   &  {-0.12} &  0      &   0 \\
\!{-0.12} &  -1.18 &  0      & 0 \\     
0.36   &    0.99 &  0.70 & 0.04 \\
\!-2.19 &  -0.04 &  0.04 & 0.17
\end{smallmatrix*}\hspace{-0.1ex}\biggr]
\!,\!  
\biggl[\!\begin{smallmatrix*}[r]
         0  & -0.97\\
         0  & -0.78\\
         0  & -0.44\\
   \!  -1.12 &   0.24\\
\end{smallmatrix*}\hspace{-0.1ex}\biggr]
\!,\!  
\biggl[\!\begin{smallmatrix*}[r]
         0     &    0 \\
         0     &    0 \\
    \! 1.57     &    0 \\
         0     & -0.49
\end{smallmatrix*}\hspace{-0.1ex}\biggr]
\\
(A_1,B_1) &= \biggl[\begin{smallmatrix*}[r] 
    -3 &    1 &    0 &    0  \\
    -2 &   -9 &    0 &    0 \\   
    -8 &    6 &   -6 &   -8 \\  
    -3 &    1 &   -8 &   -8  
\end{smallmatrix*}\biggr] \times 10^{-2}
,
\biggl[\begin{smallmatrix*}[r] 
 -4 &   -7\\
3 &   -7\\
 -4 &   -9\\
  5 &    8
\end{smallmatrix*}\biggr] \times 10^{-2}
\\
(A_2,B_2) &= \biggl[\begin{smallmatrix*}[r] 
    -7 &   -4 &    0 &    0 \\
     7 &    1 &    0 &    0 \\
     5 &   -8 &  -10 &   -7 \\
    -4 &    0 &    2 &    1 
\end{smallmatrix*}\biggr] \times 10^{-2}
,
\biggl[\begin{smallmatrix*}[r] 
   5 &    8\\
   -5 &    7\\
   -7 &    0\\
  8 &    2
\end{smallmatrix*}\biggr] \times 10^{-2}
\\
(A_3, B_3) &= \biggl[\begin{smallmatrix*}[r] 
     5 &    6 &    0 &    0 \\
    -8 &   10 &    0 &    0 \\
     0 &   -5 &   -3 &   -7 \\
    -6 &    8 &    8 &    1 
\end{smallmatrix*}\biggr] \times 10^{-2}
\biggl[\begin{smallmatrix*}[r] 
 -7 &   -6\\
   7 &   -5\\
 -3 &   -2\\
 10 &    5
\end{smallmatrix*}\biggr] \times 10^{-2}
\end{align*}
and true parameter vector $\smash{\theta^* = \begin{bmatrix} {-0.5} & {-0.152} &  0.44 \end{bmatrix}^\top}$. 
The system constraints are $\| x_t \|_\infty \leq 1$, $\|u_t\|_\infty \leq 1$, the cost matrices are $Q=\Id_{4\times 4}$, $R = \Id_{2\times 2}$,
and the disturbance sequence $\{w_0,w_1,\ldots\}$ is uniformly distributed on ${\W} = \{ w  \in \mathbb{R}^2 :  w^\top P_w w \leq 1\}$ with $P_w = \bigl[\begin{smallmatrix} 2.18 & 0.07 \\ 0.07 & 2.19\end{smallmatrix}\bigr]\times10^4$. 

This system satisfies Assumptions~\ref{ass:disturbance} and \ref{ass:param_uncertainty}, but the structure of $F$ and the block-lower-triangular structure of $A_0$, $A_1$, $A_2$ and $A_3$ 
ensures that the system is not reachable from the disturbance input $w_t$. Therefore persistence of excitation cannot be achieved through the action of the random disturbance $w_t$ alone, whereas $u_t$ may be made persistently exciting with a suitably chosen control law.
%

A robust MPC law was used with ellipsoidal tubes bounding predicted states. 
The terminal set is defined as a robustly invariant ellipsoidal set $\X_T=\{x: x^\top P x \leq 1\}$, and the tube cross-sections $\X_{k|t} = \{x : \smash{(x-z_{k|t})^\top} P (x- z_{k|t}) \leq \beta_{k|t}\}$ are defined using the same matrix $P$, with $\{\beta_{0|t},\ldots,\beta_{N|t}\}$ and $\{z_{0|t},\ldots,z_{N|t}\}$ as decision variables in problem $\P$, where 
\[
z_{k+1|t} = A_K(\bar{\theta}_t) z_k + B(\bar{\theta}_t) v_{k|t}, \ k\in\ZZ_{[0,N-1]}, 
\]
and $\bar{\theta}_t = \arg\min_{\bar{\theta}} \max_{\theta\in \Theta_t } \| \bar{\theta} - \theta\|$.
Simulations were performed in Matlab on a 2.9\,GHz Intel Core i7 processor with 16\,GB RAM. The offline MPC design steps and the online optimisation in step (c) of Algorithm~\ref{alg:sampling} were formulated using Yalmip~\citep{Yalmip} and solved using Mosek~\citep{Mosek}. Identical sequences of random disturbances were used in parallel tests of different control strategies.
By omitting step (d) from Algorithm~\ref{alg:sampling} we obtain an adaptive robust MPC law without PE guarantees, and Algorithm~\ref{alg:sampling} reduces to a robust MPC law with a fixed parameter set if steps (b) and (d) are omitted.

\begin{figure}[h]
    \centering
    \includegraphics[scale=0.5, trim = 7mm 1mm 13mm 8mm, clip=true]{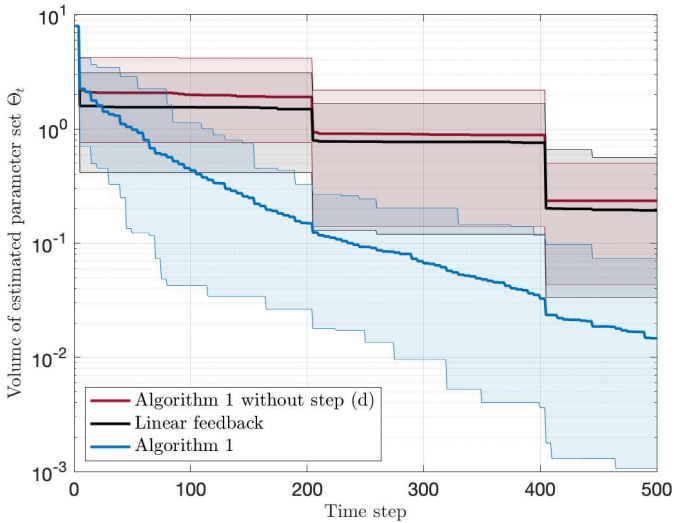}
    \caption{Parameter set volume $\vol(\Theta_t)$ for Algorithm~\ref{alg:sampling} with and without step (d), and linear feedback $u_t=Kx_t$.}
    \label{fig:volume_sampling}
\end{figure}

Figure \ref{fig:volume_sampling} shows the volume of the estimated parameter set with Algorithm~\ref{alg:sampling} (with and without step (d)) and with linear feedback. 
\bluetext{The model state switches to new initial conditions at $t = 200$ and $t = 400$ in order to illustrate the effect of transients, which could arise, for example, after reference changes in tracking problems.}
Each simulation was repeated with 30 disturbance sequences; mean parameter set volumes are shown by solid lines and ranges are shown by shaded regions. 
In agreement with Lemma \ref{lem:convergence} and Theorem \ref{thm:pe_closed_loop}, Alg.~\ref{alg:sampling} (blue) results in parameter convergence, while the adaptive MPC law of Alg.~\ref{alg:sampling} without step (d) (red) and linear feedback (black) are unable to improve their parameter estimates except during the transient responses around $t = 200$ and $t=400$.

\begin{figure}[h]
    \centering
    \includegraphics[scale=0.5, trim = 5mm 0mm 14mm 10mm, clip=true]{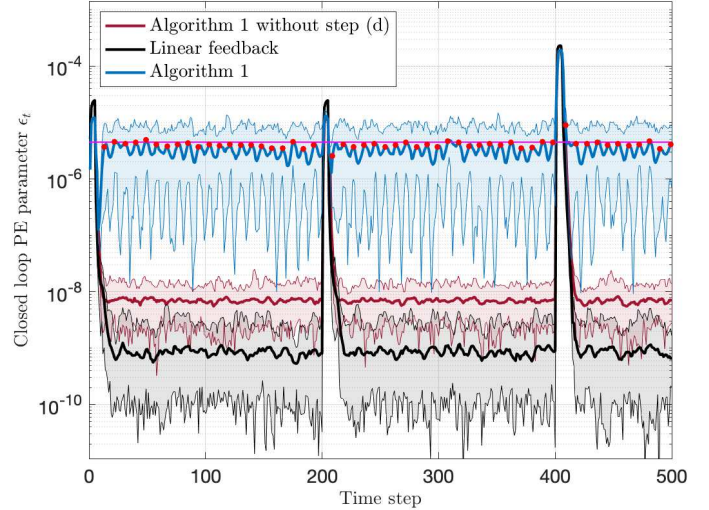}
    \caption{Closed loop PE coefficients $\epsilon_{t}$ for Algorithm~\ref{alg:sampling} with and without step (d), and linear feedback $u_t=Kx_t$.}
\label{fig:cl_PE_coeff}
\end{figure}

The convergence rates of the estimated parameter sets  in Figure~\ref{fig:volume_sampling} are explained by the measures of persistence of excitation shown in Figure~\ref{fig:cl_PE_coeff}. Here the closed loop PE coefficient at time $t$ is the largest scalar $\epsilon_t$ satisfying 
\[
\sum_{k=t-N_u+1}^{t}\Phi_k^\top \Phi_{k} \succeq \epsilon_t\Id .
\]
The solid lines show the values of $\epsilon_t$ averaged over 30 uncertainty realisations, and the shaded regions show the ranges of values observed for $\epsilon_t$. 
With the exception of short periods following the times when the system switches to a new initial condition, $\epsilon_t$ is typically orders of magnitude greater when step (d) is included in Algorithm~\ref{alg:sampling}, and $\epsilon_t$ for linear feedback is effectively zero (within limits on numerical precision).
The threshold employed in step (d) was $\epsilon = \smash{10^{-6}}$, the number samples used to compute $\delta_{\kappa|t}$ and $\smash{\hat{\delta}_{\kappa|t}}$ in step (d) was $N_s = 20$, the PE window length was $N_u = 5$, and the MPC prediction horizon was $N=10$.  

The pink line  in Figure~\ref{fig:cl_PE_coeff} shows the value of $\epsilon_\Phi$, which was defined in (\ref{eq:epsilon_Phi}) 
as a lower bound on the minimum eigenvalue of $\EE_{N|t} [ \Psi_{N-1|t}(\hat{\bf{v}}_{t},\zeta) ]$ for all $\zeta \in\Theta_t\times\smash{\W^N}$. For convenience we define $\epsilon_\Phi$ here by minimizing over $\Theta_0$, considering the contribution from $s_{k|t}$, $k\geq N$, and neglecting the contribution from the free response starting from $x_{N|t}$.
The sequence $\{s_{0|t},s_{1|t},\ldots\}$ is uniformly distributed on ${\S} = \{ s :  \smash{s^\top P_s s} \leq 1\}$ with $P_s = \bigl[\begin{smallmatrix} 1.45 & 0.16 \\ 0.16 & 1.3 \end{smallmatrix}\bigr]\times10^4$. 
The values of $\epsilon_t$ for Algorithm~\ref{alg:sampling} at times when 
$\kappa = {-N_u+1}$ are shown by red circles in Figure~\ref{fig:cl_PE_coeff}. 
Clearly the empirical expected value of $\epsilon_t$ at these times has a non-zero probability of exceeding $\epsilon_\Phi$,
in agreement with  Theorem~\ref{thm:pe_closed_loop}. 


\begin{figure}[h]
    \centering
    \includegraphics[scale=0.5, trim = 8mm 3mm 15mm 9mm, clip=true]{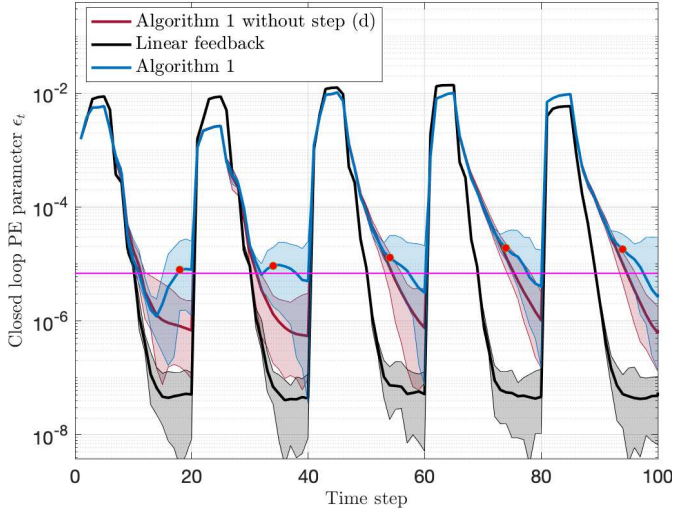}
    \caption{Closed loop PE coefficients $\epsilon_{t}$ for Algorithm~\ref{alg:sampling} with and without step (d), and linear feedback $u_t=Kx_t$.}
\label{fig:cl_PE_coeff2}
\end{figure}

To investigate the effectiveness of the PE check in Algorithm~\ref{alg:sampling} during transient responses,
we now consider switching between randomly chosen initial conditions at intervals of 20 time steps. 
Figure~\ref{fig:cl_PE_coeff2} shows that Algorithm~\ref{alg:sampling} with step (d) (in blue) ensures that the PE condition (\ref{eq:pe_condition}) is met more reliably than either Algorithm~\ref{alg:sampling} without step (d) (red) or linear feedback (black). Here the PE check in step (d) failed at 32\% of time steps on average (reduced from 84\% in Figure~\ref{fig:volume_sampling}) and linear feedback violated input or state constraints at 25\% of time steps.

\begin{figure}[h]
    \centering
    \includegraphics[scale=0.5, trim = 10mm 3mm 15mm 10mm, clip=true]{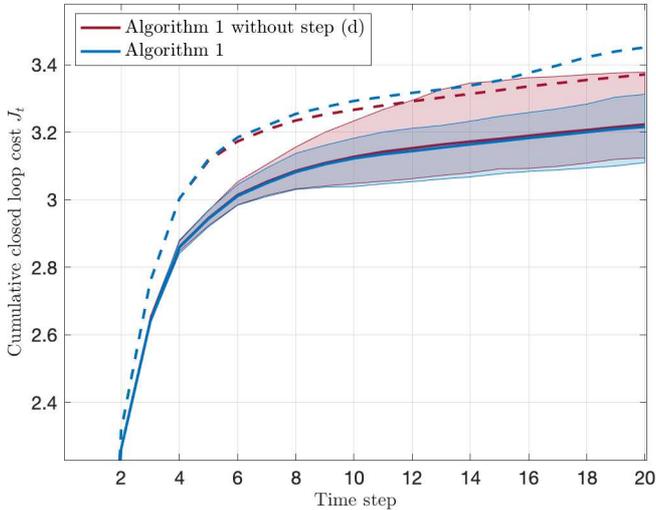}
    \caption{Closed loop cumulative cost $J_{t}$ for Algorithm~\ref{alg:sampling} with and without step (d). Dashed lines: during parameter set estimation, Solid lines: constant parameter set.}
\label{fig:cumulative_cost}
\end{figure}

To compare performance, we consider the cumulative cost at time $t$, $J_t  \defeq \sum_{k = 0}^{t} (\| x_k\|_Q^2+ \| u_k\|_R^2)$. Figure~\ref{fig:cumulative_cost} shows the evolution of $J_t$ in the first 20 times steps of Figure~\ref{fig:cl_PE_coeff2} while the parameter set is updated (dashed lines), and when $\Theta_t$ is frozen at the value of $\Theta_{200}$ (solid lines).
Comparing the blue and red dashed lines, the cumulative cost rises slightly faster when step (d) is included, indicating that the guarantee of parameter convergence provided by step (d) is obtained at the expense of slightly worse performance. However, the solid lines show that the more accurate parameter estimates obtained using Algorithm~\ref{alg:sampling} with step (d)  provide improved cumulative costs since the range of costs (shaded regions in Fig.~\ref{fig:cumulative_cost}), is smaller and the worst case cost is reduced by 5\%. 

\begin{table}[h]
\caption{Mean \{max\} computation times (ms)}
\label{tab:computation}
\vspace{-2mm}
\centerline{\begin{tabular}{@{} l | r r r @{}}
\hline
Model sizes and Alg.\,\ref{alg:sampling} parameters &  step (b)  & step (c) & step (d)\\
\hline
\begin{tabular}{@{}c@{}}
$(n_x,n_u,p) = (4,2,3)$ \\
$(N,N_u,N_s) = (10,5,20)$
\end{tabular}
& 
\begin{tabular}{@{}r@{}}
12.0\\ 
\{14.3\}
\end{tabular}
& 
\begin{tabular}{@{}r@{}}
12.9\\ 
\{17.2\}
\end{tabular}
& \begin{tabular}{@{}r@{}}
4.7\\ 
\{8.9\}
\end{tabular}
\\
\hline
\begin{tabular}{@{}c@{}}
$(n_x,n_u,p) = (4,2,{\bf 5})$ \\
$(N,N_u,N_s) = (10,5,20)$
\end{tabular}
& \begin{tabular}{@{}r@{}}
22.0\\ 
\{25.1\}
\end{tabular}
& \begin{tabular}{@{}r@{}}
43.2\\ 
\{51.8\}
\end{tabular}
& \begin{tabular}{@{}r@{}}
17.7\\ 
\{33.0\}
\end{tabular}
\\
\hline
\begin{tabular}{@{}c@{}}
$(n_x,n_u,p) = ({\bf 8},{\bf 3},5)$ \\
$(N,N_u,N_s) = ({\bf 20},{\bf 10},20)$
\end{tabular}
& \begin{tabular}{@{}r@{}}
26.4\\ 
\{33.8\}
\end{tabular}
& \begin{tabular}{@{}r@{}}
89.5\\ 
\{107.0\}
\end{tabular}
& \begin{tabular}{@{}r@{}}
23.1\\ 
\{42.5\}
\end{tabular}
\\
\hline
\end{tabular}}
\end{table}

The dependence of computation times on plant dimensions and algorithm parameters for a randomly chosen set of plant models is shown in Table~\ref{tab:computation}.

\section{Conclusions}
An adaptive robust MPC algorithm is proposed, providing a persistently exciting control law while robustly satisfying input and state constraints. Persistence of excitation is ensured by a sample-based check on the solution of the MPC optimisation problem, and by using a fallback control law if this check fails. The approach provides is computationally tractable and robustly stable.


%


\bibliography{sample1.bib}

\end{document}